\documentclass[11pt,fleqn]{article}
\usepackage{amsfonts}
\newcommand{\ol}{\setlength{\itemsep}{0pt.}\begin{enumerate}}
\newcommand{\eol}{\end{enumerate}\setlength{\itemsep}{-\parsep}}
\newcommand{\ignore}[1]{}
\title{Edge-isoperimetric inequalities and influences}
\author{Dvir Falik\thanks{School of Computer Science and
    Engineering, Hebrew University, Jerusalem, Israel.}
\and Alex Samorodnitsky\thanks{School of Computer Science and
    Engineering, Hebrew University, Jerusalem, Israel.}}

\begin{document}
\date{}
\maketitle
 
 
\newtheorem{THEOREM}{Theorem}[section]
\newenvironment{theorem}{\begin{THEOREM} \hspace{-.85em} {\bf :} 
}%
                        {\end{THEOREM}}
\newtheorem{LEMMA}[THEOREM]{Lemma}
\newenvironment{lemma}{\begin{LEMMA} \hspace{-.85em} {\bf :} }%
                      {\end{LEMMA}}
\newtheorem{COROLLARY}[THEOREM]{Corollary}
\newenvironment{corollary}{\begin{COROLLARY} \hspace{-.85em} {\bf 
:} }%
                          {\end{COROLLARY}}
\newtheorem{PROPOSITION}[THEOREM]{Proposition}
\newenvironment{proposition}{\begin{PROPOSITION} \hspace{-.85em} 
{\bf :} }%
                            {\end{PROPOSITION}}
\newtheorem{DEFINITION}[THEOREM]{Definition}
\newenvironment{definition}{\begin{DEFINITION} \hspace{-.85em} {\bf 
:} \rm}%
                            {\end{DEFINITION}}
\newtheorem{EXAMPLE}[THEOREM]{Example}
\newenvironment{example}{\begin{EXAMPLE} \hspace{-.85em} {\bf :} 
\rm}%
                            {\end{EXAMPLE}}
\newtheorem{CONJECTURE}[THEOREM]{Conjecture}
\newenvironment{conjecture}{\begin{CONJECTURE} \hspace{-.85em} 
{\bf :} \rm}%
                            {\end{CONJECTURE}}
\newtheorem{MAINCONJECTURE}[THEOREM]{Main Conjecture}
\newenvironment{mainconjecture}{\begin{MAINCONJECTURE} \hspace{-.85em} 
{\bf :} \rm}%
                            {\end{MAINCONJECTURE}}
\newtheorem{PROBLEM}[THEOREM]{Problem}
\newenvironment{problem}{\begin{PROBLEM} \hspace{-.85em} {\bf :} 
\rm}%
                            {\end{PROBLEM}}
\newtheorem{QUESTION}[THEOREM]{Question}
\newenvironment{question}{\begin{QUESTION} \hspace{-.85em} {\bf :} 
\rm}%
                            {\end{QUESTION}}
\newtheorem{REMARK}[THEOREM]{Remark}
\newenvironment{remark}{\begin{REMARK} \hspace{-.85em} {\bf :} 
\rm}%
                            {\end{REMARK}}
 
\newcommand{\thm}{\begin{theorem}}
\newcommand{\lem}{\begin{lemma}}
\newcommand{\pro}{\begin{proposition}}
\newcommand{\dfn}{\begin{definition}}
\newcommand{\rem}{\begin{remark}}
\newcommand{\xam}{\begin{example}}
\newcommand{\cnj}{\begin{conjecture}}
\newcommand{\mcnj}{\begin{mainconjecture}}
\newcommand{\prb}{\begin{problem}}
\newcommand{\que}{\begin{question}}
\newcommand{\cor}{\begin{corollary}}
\newcommand{\prf}{\noindent{\bf Proof:} }
\newcommand{\ethm}{\end{theorem}}
\newcommand{\elem}{\end{lemma}}
\newcommand{\epro}{\end{proposition}}
\newcommand{\edfn}{\bbox\end{definition}}
\newcommand{\erem}{\bbox\end{remark}}
\newcommand{\exam}{\bbox\end{example}}
\newcommand{\ecnj}{\bbox\end{conjecture}}
\newcommand{\emcnj}{\bbox\end{mainconjecture}}
\newcommand{\eprb}{\bbox\end{problem}}
\newcommand{\eque}{\bbox\end{question}}
\newcommand{\ecor}{\end{corollary}}
\newcommand{\eprf}{\bbox}
\newcommand{\beqn}{\begin{equation}}
\newcommand{\eeqn}{\end{equation}}
\newcommand{\wbox}{\mbox{$\sqcap$\llap{$\sqcup$}}}
\newcommand{\bbox}{\vrule height7pt width4pt depth1pt}
\newcommand{\qed}{\bbox}
\def\sup{^}

\def\H{\{0,1\}}

\def\S{S(n,w)}

\def \E{\mathbb E}
\def \R{\mathbb R}
\def \Z{\mathbb Z}

\def\<{\left<}
\def\>{\right>}
\def \({\left(}
\def \){\right)}
\def \e{\epsilon}
\def \l{\lfloor}
\def \r{\rfloor}
\def \B{{\cal B}}
\def \F{{\cal F}}
\def \G{{\cal G}}

\def \L{{\bigtriangledown}}

\def \a_i{\E d^2_{i,0}}
\def \b_i{\E d^2_{i,1}}

\def \f{f_{T|i}}
\def \fp{f_{T_i|\pi(i+1)}}

\newcommand{\rarrow}{\rightarrow}
\newcommand{\lrarrow}{\leftrightarrow}

\begin{abstract}
We give a combinatorial proof of the result of Kahn, Kalai, and Linial
\cite{kkl}, which states that every balanced 
boolean function on the $n$-dimensional boolean cube has a variable
with influence of at least $\Omega\(\frac{\log n}{n}\)$.

The methods of the proof are then used to recover additional
isoperimetric results for the cube, with improved constants.

We also state some conjectures about optimal constants and discuss
their possible implications.

\end{abstract}

\section{Introduction}
\label{Introduction}
This paper deals with isoperimetric problems on graphs. Given a graph
$G = (V,E)$, the vertex boundary of a subset $S\subseteq V$ contains
the vertices of $S$ which have neigbours outside $S$
$$
\B_v(S) = \left\{x \in S:~\exists y \in S^c \mbox{ such that} ~(x,y)
\in E \right\}
$$
The edge boundary of $S$ is the set of edges crossing from $S$ to its
complement. 
$$
\B_e(S) = \left\{(x,y)\in E:~x \in S \mbox{ and }  y \in S^c \right\}
$$
The question about the smallest possible boundary a
set of given cardinality can have is an important combinatorial
question with obvious connections to the classical isoperimetry.  
Good estimates of the minimal boundary size are also very useful in
applications. We briefly mention two. Lower bounds on the 
vertex boundary show how fast a
neighbourhood of a set has to grow when the allowed distance from the
set increases, and this leads to concentration of measure results
for Lipschitz functions on the graph \cite{milman:schechtman}. Lower 
bounds on the edge boundary suggest that a simple
random walk on the graph doesl not remain in any subset for too long, and
this leads to upper bounds on its mixing time \cite{jerrum03}. 

Early isoperimetric results on graphs include isoperimetric theorems
for the boolean cube $\H^n$. This is a graph with $2^n$ vertices
indexed by boolean strings of length $n$. Two vertices are
connected by an edge if they differ only in one coordinate. The 
metric defined by this graph is called the Hamming distance.  Two vertices $x$
and $y$ are at distance $d$ if they differ in $d$ coordinates. A
Hamming ball is a ball in this metric. A subcube is a
subset of the vertices obtained by fixing the value in some of the
coordinates. The number of fixed coordinates is called the co-dimension of
the subcube. It turns out \cite{harper}, \cite{harp, hart} that the
vertex boundary of a Hamming ball is smallest among all sets of equal
size, and the same is true for the edge boundary of a subcube. (It is
also possible to interpolate these results for all the intermediate
subset sizes.)

While many other exact vertex and edge isoperimetric results are known
(see \cite{bezrukov} for a survey), in most cases exact
results seem to be hard to obtain. In many of these cases they could
be replaced by sufficiently strong approximate
isoperimetric results \cite{js, moo, product-measures}. 

Given a solution of an isoperimetric problem, one can ask about its
stability. Namely, should a set whose boundary is not much larger than
minimal be close to the optimal set? 
Such results turn out to be especially useful and 
interesting \cite{friegut-bourgain, friedgut, kahn-kalai}. 

In this paper we focus on edge-isoperimetric questions in the boolean
cube. A major result in this area was 
obtained by Kahn, Kalai, and Linial \cite{kkl} who showed that any
balanced boolean function has a variable with large {\it
influence}. We proceed to describe this result, starting with some
background. 

For a subset $A \subseteq \H^n$ and an index $1\le i\le n$ let $I_i(A)$
be the fraction of edges in direction $i$ 
between $A$ and its complement $A^c$. This means that $2^{n-1} \cdot
I_i(A)$ counts 
the edges with one vertex in $A$ and another in $A^c$, the vertices
disagreeing in $i$-th coordinate. $\sum_{i=1}^n I_i(A)$ is
the total (normalized) cardinality of the edge boundary of $A$.

The familiar edge-isoperimetric inequality in the cube states that for
any subset $A \subseteq \H^n$ of cardinality at most $2^{n-1}$ holds
\footnote{We use natural logarithms throughout the paper.}
$$
\sum_{i=1}^n I_i(A) \ge \frac{2}{\log 2} \cdot \frac{|A|}{2^n} \log
\frac{2^n}{|A|}. 
$$
This is tight if $A$ is a subcube of (arbitrary) co-dimension $1\le t \le n$.

Let $f$ be the characteristic function of $A$, with expectation $\mu =
\E_{x \in \H^n} f(x) = \frac{|A|}{2^n}$. The edge-isoperimetric
inequality asserts that for $\mu \le 1/2$
\beqn
\label{set_isop}
\sum_{i=1}^n I_i(f) \ge \frac{2}{\log 2} \cdot \mu
\log\frac{1}{\mu}.
\eeqn
Here $I_i(f)$ stands for the influence of the $i$-th variable on
the support of $f$.\footnote{We interchange freely between a set and
its characteristic function. Whenever this does not cause confusion
we do not mention either, and simply write $I_i$.}

\ignore{
Another variation
valid for all $0 \le \mu \le 1$ is
\beqn
\label{sym_set_isop}
\E_x \sum_{y \sim x} (f(x) - f(y))^2 \ge \frac{2}{\log 2} \cdot
\mu(1-\mu) \log\frac{1}{\mu(1-\mu)}
\eeqn
This is somewhat weaker than (\ref{set_isop}) but the
isoperimetric constant 
$C_{\ref{sym_set_isop}} = \frac{2}{\log 2}$ is still tight, for a subcube
of dimension $n-1$. A benefit of this version lies in its consistency
with other isoperimetric inequalities we deal with here, which
are stated in terms of $\mu(1-\mu)$.
}
The inequality (\ref{set_isop}) has several easy proofs \cite{harp,
  hart}. The 
one most relevant to this 
discussion is by induction on dimension. To illustrate its outlay and
its simplicity, here it is (a sketch): the
base $n = 1$ is easy . Assume for dimension $n-1$ and consider the
case of dimension $n$. Write $A = A_0 \cup A_1$, where $A_i$ contains
all the elements of $A$ with $i$ in the $n$'th coordinate. Think about
$A_i$ as subsets of $(n-1)$-dimensional cube, and observe $I_n(A)
\ge \(1/2^{n-1}\) \cdot \Big | |A_0| - |A_1| \Big |$. Taking $a_i =
|A_i|$ it remains to 
check that for any nonnegative $a_0$, $a_1$ holds
$
2 \log 2 a_0 \log \frac{2^{n-1}}{a_0} + 2 \log 2  a_1 \log
\frac{2^{n-1}}{a_1} + |a_0 - a_1| \ge 2 \log 2  (a_0+a_1) \log
\frac{2^{n}}{a_0+a_1}, 
$
which is easily verified, using the properties of the logarithm.
 
Things become more complicated when we ask for more detailed
information. Interpreting the set $A$ as the set of positive outcomes
of a game with $n$ players,
the number $I_i(A)$ acquires a game-theoretic interpretation as the
{\it influence} of $i$-th player on the outcome of the game, namely
the probability that the outcome of the game remains uncertain if the
decisions of other players are chosen at random. Motivated by questions
from computational game theory Ben-Or and Linial \cite{bl} conjectured
that for any balanced game (namely $|A| = 2^{n-1}$) there is a player
with influence of at least $\Omega\(\frac{\log n}{n}\)$. 

This conjecture was proved by Kahn, Kalai, and Linial in \cite{kkl}. 
\thm
\label{kkl}
Let $f~:~\H^n \rarrow \{0,1\}$ be a boolean function with expectation
$\E f = \mu$. Then
\beqn
\label{ineq_kkl}
\sum_{i=1}^n I^2_i(f) \ge \Omega\(
\frac{\mu^2(1-\mu)^2\log^2 n}{n}\) 
\eeqn
\ethm
In particular, there is $i$ with $I_i \ge \Omega\(
\frac{\mu(1-\mu)\log n}{n}\)$. 

\cite{kkl} 
is one of the first papers to use Fourier analysis on $\Z^n_2$ in a
combinatorial setting. Rather surprisingly, a crucial tool in the 
proof is an inequality \cite{Beck, Bon,
  Gr} which is easiest to describe in Fourier-analytic terms. Let
$\left\{w_S\right\}$ be the Walsh-Fourier basis of the vector space
of real-valued functions on the cube. For a function 
$f :~\H^n \rarrow \R$, $f = \sum_{S \in \H^n} \hat{f}(S) w_s$, and a
nonnegative real $\e$,  let $T_{\e}(f) = \sum_{S \in \H^n} \e^{|S|}
\hat{f}(S) w_s$. Then
$$
\|T_{\e} f\|_2 \le \|f\|_{1 + \e^2}
$$
Following its application in \cite{kkl}, this inequality, known
(for historical reasons) as the Bonami-Beckner 
inequality, became a very important tool in combinatorics and theory
of computer science. Still it is very different from the familiar
combinatorial tools, and its appearance in the proof is somewhat
mysterious. 
Thus it seemed of interest to look for a
combinatorial proof of theorem~\ref{kkl}, possibly along the lines of
the forementioned proof of (\ref{set_isop}).
Let us mention two papers dealing with this problem along very
different routes.
The first of these papers \cite{friedgut-rodl} gives a
combinatorial (entropic) proof of the Bonami-Beckner 
inequality for $\e = \sqrt{3}/3$. This special case is already strong
enough to be instrumental in the
proof of theorem~\ref{kkl}. The second paper
\cite{talagrand-boundar-infl} presents an inductive  
proof that the maximal influence of a balanced function is at least
$\Omega\(\frac{\log^{\alpha}(n)}{n}\)$ for some $0 < \alpha < 1$.  

In this paper we give a fully combinatorial proof of
theorem~\ref{kkl}. After completing our work,
we learned that a very similar proof was recently obtained by
\cite{rossignol}. 

We start with a functional
form of inequality (\ref{set_isop}). For a nonnegative function $f:~\H^n
\rarrow \R$ 
\beqn
\label{func_isop}
\E_x \sum_{y \sim x} (f(x) - f(y))^2 \ge 2 \cdot \E f^2 \log \frac{\E
  f^2}{\E^2 f}
\eeqn
Functional forms of isoperimetric inequalities are widely used in
local theory of Banach spaces \cite{latala}. They turn
out to be useful in our setting too. We show theorem~\ref{kkl} to be a
simple consequence of inequality (\ref{func_isop}). 

This inequality can be proved by induction on dimension (see
  Appendix A), similarly to
(\ref{set_isop}), though the proof is somewhat more complicated. The
isoperimetric constant $C_{\ref{func_isop}} = 2$ is tight, if we want
it to be independent of the dimension. In
section~\ref{Construction} we give examples of functions satisfying
(\ref{func_isop}) with equality if the constant $2$ is replaced by $2
  + o_n(1)$. These are symmetric functions (a function $f$ on the cube
  is symmetric if $f(x)$ depends only on the distance of
  $x$ from zero) closely related to a classical family of orthogonal
  polynomials of discrete variable - the Krawchouk polynomials.

We also suggest a reason behind the relevance of the Bonami-Beckner
inequality. It is well-known that this inequality is equivalent to the
logarithmic Sobolev inequality 
\beqn
\label{ineq_Sob}
\E_x \sum_{y \sim x} (f(x) - f(y))^2 \ge 2 \cdot Ent\(f^2\) = 2\cdot \(\E f^2
\log f^2 - \E f^2 \log \E f^2\)
\eeqn 
in the cube. It turns out that this inequality implies inequality
(\ref{func_isop}). 
In this sense Bonami-Beckner's inequality can be thought of as a refined form
of the edge-isoperimetric inequality in the discrete cube. 

The actual result we prove seems to be somewhat stronger than
theorem~\ref{kkl}. We show that for a boolean function with
expectation $\mu$ holds
\thm
\label{boolean}
\beqn
\label{ineq_main_boolean}
\sum_{i=1}^n I^2_i(f) \ge 4\mu(1-\mu) \exp\left\{-\frac{1/2}{\mu(1-\mu)}
\sum_{i=1}^n I_i(f) \right\}    
\eeqn
\ethm
This inequality implies (\ref{ineq_kkl})
with a constant $c_{\ref{ineq_kkl}} = 4$, recovering the estimate of
\cite{kkl}. 

Theorem~\ref{boolean} is a special case of our main result,
theorem~\ref{main}, stated and proved in
section~\ref{Main}. This theorem presents a more general inequality
valid for real-valued functions on the discrete cube endowed with an
arbitrary measure. The theorem and the
approach used in its proof seem to provide convenient tools for
dealing with a certain type of isoperimetric statements in the
cube. We illustrate this by giving simple proofs of two results from
\cite{friedgut} and \cite{friedgut-kalai}, with better
isoperimetric constants. 

It should be mentioned that inequality (\ref{ineq_main_boolean}), in
its turn, is implied, {\it up to a
constant in the exponent}, by an inequality
of Talagrand \cite{tal}. A special case of this inequality asserts
that for a boolean function $f$ with expectation $\mu$ 
\beqn
\label{ineq_talagrand}
\sum_{i=1}^n \frac{I_i}{\log\(e/I_i\)} \ge \Omega\(
\mu (1-\mu)\)
\eeqn
It is not hard to see that this gives (\ref{ineq_main_boolean}) if
$1/2$ in the exponent is replaced by a sufficiently large constant.

Next, we focus our attention on the best possible constants for the
above-mentioned inequalities. Specifically, we are interested in
the exact constant
$C_{\ref{ineq_main_boolean}}$  that should appear in the exponent in
the right hand side of (\ref{ineq_main_boolean}). We point out an
interesting phenomenon in that obtaining the (conjectured) optimal
constant for this inequality would lead to a stability result for the
basic inequality (\ref{set_isop}). 

To be more specific, 
by theorem~\ref{boolean} $C_{\ref{ineq_main_boolean}} \le \frac12$. 
On the other hand, taking $f$ to be a characteristic function of a
subcube of large co-dimension, shows $C_{\ref{ineq_main_boolean}} \ge
\frac{\log 2}{2}$. We believe the lower bound to be the right one.
\cnj
\label{main_cnj}
$$
C_{\ref{ineq_main_boolean}} = \frac{\log 2}{2}
$$
\ecnj
In particular, we conjecture small subcubes to be (nearly)-isoperimetric
sets for this inequality.

If conjecture~\ref{main_cnj} holds, this would, in particular, give
the optimal constant $C_{\ref{ineq_kkl}} = \frac{4}{\log^2 2} \approx
8.3$ in inequality (\ref{ineq_kkl}). It is easy to see that 
$C_{\ref{ineq_kkl}} \ge
  \frac{1}{C^2_{\ref{ineq_main_boolean}}}$.\footnote{More precisely,
  $C_{\ref{ineq_kkl}} \ge 
  \frac{1}{C^2_{\ref{ineq_main_boolean}}} - o_n(1)$. Here and in the
  rest of this paper we ignore negligible factors when comparing constants.}
Therefore the conjecture would imply $C_{\ref{ineq_kkl}} \ge
\frac{4}{\log^2 2}$. It turns out that the best known
candidate to be an isoperimetric function for inequality (\ref{ineq_kkl}), the
``tribes'' function of Ben-Or and Linial \cite{bl} indeed shows
$C_{\ref{ineq_kkl}} \le \frac{4}{\log^2 2}$.
\footnote{Choosing the tribe size appropriately, so that the
  expectation is small.} We are grateful to Amites
Sarkar \cite{asarkar} for pointing this out to us.

Now, consider functions which are nearly
isoperimetric in the sense of the basic inequality
(\ref{set_isop}). Kahn and Kalai conjecture (\cite{kahn-kalai}) that such
functions behave similarly to subcubes, in the following precise sense.
\cnj \cite{kahn-kalai}
Let $K > 0$ be a real number. There are positive real numbers $K',
\delta$ depending on $K$ such that the following assertion holds: 
If a monotone boolean function $f:~\H^n \rarrow \{0,1\}$ with expectation $\mu
\le 1/2$ satisfies
\beqn
\label{cnj_kk}
\sum_{i=1}^n I_i \le K \cdot \frac{2}{\log 2} \cdot \mu \log
\frac{1}{\mu}
\eeqn
then there is a set of at most $K' \cdot \log \frac{1}{\mu}$
coordinates such that the expectation of $f$ restricted to the subcube 
obtained by setting all the coordinates in this set to $1$ is at least
$(1 + \delta) \cdot \mu$.  
\ecnj
It seems that a weaker version of this conjecture, which claims the
same conclusion from a stronger assumption that the multiplicative
factor $K$ in (\ref{cnj_kk}) is close to $1$, i.e., $K = 1 + \e$ for a
small $\e > 0$, is also interesting \cite{kalai-personal}. 

Here we prove an even weaker result in this direction, conditioned on
conjecture~\ref{main_cnj}. Let
$o_{\mu,\e}(1)$ denote a quantity which goes to zero when {\it both}
$\mu$ and $\e$ do. 
\pro
\label{cond_kk}
Assume conjecture~\ref{main_cnj}. 
Let $f:~\H^n \rarrow \H$ be a monotone boolean function with expectation $\mu
\le 1/2$, and assume
$$  
\sum_{i=1}^n I_i \le (1 + \e) \cdot \frac{2}{\log 2} \cdot \mu \log
\frac{1}{\mu}
$$
Then there is a set of $O\(\(\frac{1}{\mu}\)^{\(1+o_{\mu,\e}(1)\)\cdot \e}\)$
coordinates such 
that the expectation of $f$ restricted to the subcube 
obtained by setting all the coordinates in this set to $1$ is at least
$2 \mu$. 
\epro

We conclude this section by saying a few words about a possible
approach to the proof of
conjecture~\ref{main_cnj}. We will say more about this in
section~\ref{Corollaries}. The main step in the proof of
theorem~\ref{boolean} is 
a variant of the logarithmic Sobolev inequality for the
discrete cube. This inequality applies to general real-valued
functions on the cube, and is tight with constant $c=2$. To
prove the conjecture we need to take into account the specific
structure of boolean functions. The familiar approach using
tensorization does not seem to be convenient 
for this. We give a proof of the inequality for general functions
which works by induction on the dimension, 
similar to the proof of (\ref{set_isop}), and seems to be more
conducive for this purpose.

The paper is organized as follows: in the next section we prove the
main theorem \ref{main}. In
section~\ref{Corollaries} several corollaries
are derived from theorem~\ref{main}, and the main technical
conjecture is stated. Section~\ref{Construction}  
constructs nonnegative real-valued functions which are almost
isoperimetric for inequality (\ref{ineq_main}) and hence for several other
inequalities in this paper, including (\ref{func_isop}). Inductive proofs of
(\ref{func_isop}) and a logarithmic Sobolev inequality for the cube
are given in the Appendices.

\section{The main theorem}
\label{Main}
We start with some definitions and notation.

Let $\F_j$, for $0 \le j \le n$, be the algebra of subsets of $\H^n$
generated by the first $j$ bits. More precisely, $\F_j$ is generated
by the atoms $\left\{A_{\e_1...\e_j}:~\e_i \in \H\right\}$ 
where $A_{\e_1...\e_j} = \{x:~x_1 = \e_1...x_j = \e_j\}$.
Then $\{\F_j\}$ is an increasing sequence of algebras. In particular $\F_0 =
\{\emptyset, \H^n\}$ and $\F_n = 2^{\H^n}$.

For a function $f:~\H^n \rarrow \R$, let $f_i = \E\(f|\F_i\)$, the
conditional expectation of $f$ given the algebra $\F_i$. This means
that $f_i(x)$ is the average of $f$ over the points $y$ that
coincide with $x$ in the first $i$ coordinates. In
particular, $f_0 = \E f$, $f_n = f$. The sequence $f_0,...,f_n$ is a
martingale with respect to $\{\F_j\}$.
\footnote{Essentially the only martingale property we use is the fact
 that conditional expectation is an orthogonal projection on a
 subspace.}
Let $d_i$, $i=1...n$ be the
sequence of martingale differences. $d_i = f_i - f_{i-1}$, $1\le i\le
n$. 

Let ${\cal E}(f,f) = \E_x
\sum_{y \sim x} (f(x) - f(y))^2$. Let us mention that ${\cal E}(f,g) =
E_x \sum_{y \sim x} (f(x) - f(y))(g(x) - g(y))$ 
is sometimes called the canonical Dirichlet form on $\H^n$ \cite{bobkov}.

The following lemma is simple and well-known
\cite{elchanan-personal}. For completeness, we will give a proof at
the end of this section. 
\lem
\label{energy-additive}
$$
{\cal E}(f,f) = \sum_{i=1}^n {\cal E}\(d_i,d_i\)
$$
\elem

Let $\mu$ be a measure on $\H^n$. Let $C$ be the best
constant in the logarithmic Sobolev inequality for $\H^n$ with
$\mu$. This is to say that $C$ is maximal such that for any function
$f:~\H^n \rarrow \R$ holds ${\cal E}(f,f) \ge C \cdot  Ent\(f^2\)$. 

For a function $f:~\H^n \rarrow \R$, let $\sigma^2(f) = \E f^2 -
\E^2 f$.
 
Our main result is:
\thm
\label{main}
\beqn
\label{ineq_main}
\sum_{i=1}^n \E^2 |d_i| \ge \sigma^2(f) \exp\left\{-\frac{{\cal
E}(f,f)}{C \sigma^2(f)}\right\}
\eeqn
\ethm
\prf
First a simple lemma.
\lem
\label{sob-isop}
For a nonnegative function $f$ holds $Ent(f^2) \ge \E f^2 \log
\frac{\E f^2}{\E^2 f}$.
\elem
\prf (Of the lemma)
Since both sides of the inequality are $2$-homogeneous, this amounts
to showing $\E f^2 \log f^2 + \log \E^2(f) \ge 0$, given $\E f^2 = 1$. 
This is the same as  
$\log \E f - \E f^2 \log \frac{1}{f} \ge 0$. And this is true since
logarithm is concave.
\eprf

Now we can conclude the proof of theorem~\ref{main}.
Using (\ref{ineq_Sob}) and lemma~\ref{sob-isop},
$$
{\cal E}(f,f)  = \sum_{i=1}^n {\cal E}\(d_i, d_i\) \ge C  \cdot \sum_{i=1}^n
Ent\(d^2_i\) \ge C  \cdot \sum_{i=1}^n \E d^2_i \log \frac{\E d^2_i}{\E^2
  |d_i|}.
$$
Observe $\sum_{i=1}^n \E d^2_i = \E f^2 - \E^2 f=: \sigma^2(f)$. By the
convexity of the minus logarithm, the last sum is 
$$
C \sigma^2(f) \sum_{i=1}^n \frac{\E d^2_i}{\sigma^2(f)} \log \frac{\E
  d^2_i}{\E^2 |d_i|} \ge C \sigma^2(f) \log
\(\frac{\sigma^2(f)}{\sum_{i=1}^n \E^2 |d_i|}\).
$$ 
\eprf

We remark that instead of the logarithmic Sobolev inequality
(\ref{ineq_Sob}) it is possible to use isoperimetric
inequality (\ref{func_isop}) directly. Hence the 
logarithmic Sobolev constant $C$ in the statement can be replaced by
potentially bigger isoperimetric constant $C'$. 

\noindent {\bf Proof of lemma~\ref{energy-additive}}\\
For a function $g$ and an index $1\le i \le n$ let
$g^i$ be a function defined by $g^i(x) = g(x\oplus e_i)$. Here $e_i$
is the vector with $1$ in $i$'th coordinate, and zero in the other
coordinates. Note that $g^i$ is $\F_k$-measurable iff $g$ is.
$$
{\cal E}\(d_j,d_j\) = \sum_{i=1}^n \|d_j -
d^i_j\|^2_2 = \sum_{i=1}^n \<d_j - d^i_j, d_j - d^i_j \> = 
\sum_{i=1}^n \<d_j - d^i_j, (f_j - f_{j-1}) - \(f^i_j - f^i_{j-1}\)\>
= 
$$
$$
\sum_{i=1}^n \<d_j - d^i_j, f_j - f^i_j\> = 
\sum_{i=1}^n \<\(f_j - f^i_j\) - \(f_{j-1} - f^i_{j-1}\), f_j -
f^i_j\> = 
$$
$$
\sum_{i=1}^n \<f_j - f^i_j, f_j - f^i_j\> - 
\sum_{i=1}^n \<f_j - f^i_j, f_{j-1} - f^i_{j-1}\> =
$$
$$ 
\sum_{i=1}^n \<f_j - f^i_j, f_j - f^i_j\> - 
\sum_{i=1}^n \<f_{j-1} - f^i_{j-1}, f_{j-1} - f^i_{j-1}\>\ = 
{\cal E}\(f_j,f_j\) - {\cal E}\(f_{j-1},f_{j-1}\).
$$
The proof is concluded by observing ${\cal E}\(f_0,f_0\) = 0$.
\eprf

\section{Some corollaries for product measures}
\label{Corollaries}
In this section we derive theorem~\ref{kkl},
proposition~\ref{cond_kk}, and the theorems of Friedgut and
Friedgut-Kalai from theorem~\ref{main}. We also state our main
technical conjecture (\ref{ineq_cnj_sobolev}).

We will assume the measure $\mu$ to be a product probability measure,
that is $\mu = 
\otimes_{k=1}^n \mu_k$, with $\mu_k(1) = p_k$, and $\mu_k(0) = 1 - p_k$. 
In this case $\E |d_i|$ has a simple upper bound.
\lem
\label{d_i_bound}
For a product measure $\mu$, 
$$
\E |d_i| \le 2p_i(1-p_i) \cdot \E_x \Big |f(x) - f(x\oplus e_i)\Big |
$$
\elem
\prf
$$
\E |d_i| = \E |f_i - f_{i-1} | = \E \Big | \E\(f | \F_i\) - \E\(f |
\F_{i-1}\) \Big |.
$$
Let $\G_i$ be the algebra of subsets of $\H^n$ generated by all the
bits but $j$. That is for $x$ with $x_i = 0$ holds $\E\(f | \G_i\)(x) = \E\(f |
\G_i\)(x\oplus e_i) = (1-p_i) f(x) + p_i f(x \oplus 
e_i)$.
  
Then for a product measure $\mu$ holds $\E\(\E\(f | \F_i\) \Big | \G_i\) =
\E\(f | \F_{i-1}\)$. Therefore
$$
\E \Big | \E\(f | \F_i\) - \E\(f |\F_{i-1}\) \Big | =  
\E \Big | \E \( f - \E\(f | \G_i\) \Big | \F_i\) \Big | \le 
\E \Big |f - \E\(f | \G_i\) \Big | = 2p_i(1-p_i) \cdot \E_x \Big |
f(x) - f\(x \oplus e_i\) \Big |
$$
We have used the well-known fact that conditional expectation
decreases the $\ell_1$-norm.
\eprf
\subsection{Uniform measure}
The best constant $C$ for a cube
endowed with the uniform measure is $C = 2$. Hence the theorem gives,
for a real-valued function $f$ on the discrete cube,
$
\sum_{i=1}^n \E^2 |d_i| \ge \sigma^2(f) \exp\left\{-\frac{{\cal
    E}(f,f)}{2 \sigma^2(f)}\right\}
$
This may be somewhat simplified for monotone functions, for which $\E
|d_i| = \hat{f}(\{i\})$
$$
\sum_{i=1}^n \hat{f}^2(\{i\}) \ge \sigma^2 
\exp \left\{-\frac{1}{2\sigma^2} \sum_{i=1}^n I_i(f) \right\}.
$$
Section~\ref{Construction} presents a construction of monotone functions for
which this inequality is essentially tight.

Our main concern are boolean functions. We present several easy
implications of theorem~\ref{main} and a related conjecture.

\noindent {\bf Proof of theorem~\ref{boolean}}.

For a boolean function $f$ we
have $\E_x | f(x) - f(x \oplus e_i) | = I_i$ and using the theorem
together with lemma~\ref{d_i_bound} yields 
$
\sum_{i=1}^n I^2_i(f) \ge 4\sigma^2(f) \exp
\left\{-\frac{1}{2\sigma^2(f)} \sum_{i=1}^n I_i(f) \right\}   
$
proving (\ref{ineq_main_boolean}) and theorem~\ref{boolean}. \eprf

\xam
Let $f$ be the characteristic function of a subcube of dimension
$n-t$. Then $f$ has $t$ non-zero influences of size
$2^{-t+1}$. Assume $t$ is large enough so that $\mu = 2^{-t}$ may be
replaced with $1$, and  
$\sigma^2 = \frac{2^t - 1}{4^t}$ may be replaced with $2^{-t}$. 
Then (\ref{ineq_main_boolean}) gives
$$
\frac{4t}{4^t} \ge \frac{4}{2^t} \exp\left\{-2^{t-1} \cdot
\frac{2t}{2^t}\right\} = \frac{1}{2^t}
\exp\left\{-t\right\}~~~~or~~~~~
t 2^{-t} \ge e^{-t}.
$$ 
\exam

We conjecture (conjecture~\ref{main_cnj}) that for boolean functions a
stronger inequality  
$
\sum_{i=1}^n I^2_i(f) \ge 4\sigma^2(f) \exp
\left\{-\frac{\log 2}{2\sigma^2(f)} \sum_{i=1}^n I_i(f) \right\}   
$
holds. Such an inequality would be tight by the example above.

This inequality would follow from the following version of the
logarithmic Sobolev inequality for boolean functions. 
\cnj
\label{cnj_sobolev}
For a boolean function $f$ on the discrete cube
\beqn
\label{ineq_cnj_sobolev}
{\cal E}(f,f) \ge \frac{2}{\log 2} \cdot \sum_{i=1}^n Ent\(d^2_i\)
\eeqn
\ecnj
{\bf Discussion.} The inequality holds with a constant $c = 2$ for
real-valued functions on $\H^n$. This is proved
by applying the logarithmic Sobolev inequality to functions $d_i$,
$1\le i\le n$. This is also tight, as shown by
functions constructed in section~\ref{Construction}. 
To improve the constant for boolean functions a
different approach seems to be required, one that ``remembers'' that
$\{d_i\}$ are difference functions of a martingale defined by a
boolean function. In particular the familiar tensorization approach
might not be sufficient here since 
it does not keep track of the combinatorial structure of the functions
involved. 

We give an inductive proof of the inequality ${\cal E}(f,f)
\ge 2 \cdot \sum_{i=1}^n Ent\(d^2_i\)$ in Appendix B. This
proof seems to be better suited for handling functions $f$ with a
specific structure, such as boolean functions.

\noindent {\bf Proof of theorem~\ref{kkl}}

We show
$$
\sum_{i=1}^n I^2_i(f) \ge \(4 - o_n(1)\) \cdot
\frac{\mu^2(1-\mu)^2\log^2 n}{n} 
$$
For a boolean function $f$ with expectation $\mu$ holds $\sigma^2(f) =
\mu(1-\mu)$. Let $\e(n) \gg \frac{\log\log n}{\log n}$. There are two
cases to consider. First, $\sum_{i=1}^n I_i <  
\(2-\e(n)\) \cdot \mu(1-\mu) \frac{\log n}{n}$. In this case
(\ref{ineq_main_boolean}) implies $\sum_{i=1}^n
I^2_i(f) \gg \mu(1-\mu) \frac{\log^2 n}{n}$. 
The second case is $\sum_{i=1}^n I_i \ge \(2-\e(n)\)
\cdot \mu(1-\mu)\frac{\log n}{n}$.
The Cauchy-Schwarz inequality now implies $\sum_{i=1}^n I^2_i \ge 
\(2-\e(n)\)^2 \cdot \mu^2(1-\mu)^2 \frac{\log^2 n}{n}$. \eprf

\noindent {\bf Proof of a theorem of Friedgut}
\thm \cite{friedgut}
For a boolean function $f$ and an arbitrary $\epsilon > 0$ 
there is a function $g$ depending only on $\(\sum_{i=1}^n I_i\) \cdot
\exp\left\{\frac{\sum_{i=1}^n I_i}{\(2-o_{\e}(1)\) \epsilon}\right\}$ 
coordinates\footnote{The original proof in \cite{friedgut} has a
somewhat larger estimate for the required number of variables in
the junta. This estimate has a constant $2$ instead of $1/2$ in the
exponent.}  (a junta) such that $\|f - g\|^2_2 \le \epsilon$.
\ethm
\prf
Let $K = \sum_{i=1}^n I_i$ and take $\alpha =
\exp\left\{-\frac{K}{\(2-o_{\e}(1)\)
  \epsilon}\right\}$. The error term $o_{\e}(1)$ will be chosen later.

Without loss of generality assume the influences $I_i$ to decrease
with $i$, and let $r$ be the maximal index with $I_r \ge
\alpha$. Clearly $r \le 
\frac{K}{\alpha}$ and $\sum_{i = r+1}^n I^2_i \le K \alpha$. Take $g =
\E\(f | \F_r\)$. The function $g$ 
depends only on $r$ variables. We will show $\|f - g\|^2_2 \le
\epsilon$. 

Let $h = f - g$. Take $h_i = \E\(h | \F_i\)$, $i = 1...n$, to be the
martingale defined by $h$, and let $d_i(h)$ be its difference
functions. Then $d_i(h) = d_i(f)$ for $i > r$ and $d_i(h) = 0$
otherwise. Note that $\E h = 0$ and therefore $\sigma^2(h) =
\|h\|^2_2$. By theorem~\ref{main}
$$
K = \sum_{i=1}^n I_i = {\cal E}(f,f) \ge {\cal E}(h,h) \ge 2\|h\|^2_2
\log\left(\frac{\|h\|^2_2}{\sum_{i=r+1}^n I^2_i}\right) \ge 
2\|h\|^2_2 \log\left(\frac{\|h\|^2_2}{K \alpha}\right)
$$
Recalling the definition of $\alpha$, it is now easy to choose the
error term $o_{\e}(1)$ appropriately, so that the last inequality
implies $\|h\|^2_2 \le \e$.
\eprf

\noindent {\bf Proof of proposition~\ref{cond_kk}}

We proceed similarly to the preceding proof. Let $K = \sum_{i=1}^n
I_i$ and let $\alpha =
\mu^{1 + \(1+o_{\mu,\e}(1)\)\cdot \e}$, where 
$o_{\mu,\e}(1)$ is  an error term which goes to zero when both 
$\mu$ and $\e$ do.  We will set it later.

Assume the influences to decrease, and define the index $r$ and and
functions $g$ and $h$ as above. Now we need a simple lemma
\lem
\label{close_a_bit}
If $\|h\|^2_2 \le \mu - 2\mu^2$ then 
$$
\mbox{Pr}\left\{f(x) = 1 ~\Big |~x_1 = ... = x_r = 1\right\} \ge 2 \mu
$$
\elem
\prf
For $y  = (y_1...y_r) \in \H^r$ let $K_y$ be the subcube $\left\{x \in
\H^n:~x_1 = y_1,~...~x_r = y_r\right\}$. Let $f_y$ be the function $f$
restricted to $K_y$. Let $\mu_y = \E f_y$. Alternatively $\mu_y$ is
the value of $g = \E\(f | \F_r\)$ on $K_y$. Let ${\bf 1} \in \H^r$ be
the vector of all ones. $\mu_{\bf 1}$ is the quantity we want to lower
bound. Since $f$ is a monotone function,
so is $g$. In particular $\mu_{\bf 1}$ is the largest among all
$\mu_y$. 

We have $\E_y \mu_y = \E g = \E f = \mu$. On the other hand,
$$
\mu - 2\mu^2 \ge \|h\|^2_2 = \|f - g\|^2_2 = \E_y \mu_y\(1-\mu_y\)
$$
Therefore $2 \mu^2 \le \E_y \mu^2_y \le \mu_{\bf 1} \cdot \E_y \mu_y =
\mu_{\bf 1} \cdot \mu$.
\eprf

Now, by conjecture~\ref{main_cnj}
$$
(1+ \e) \frac{2}{\log 2} \mu \log \frac{1}{\mu} \ge \sum_{i=1}^n I_i =
{\cal E}(f,f) \ge {\cal E}(h,h) \ge 
$$
$$
\frac{2}{\log 
  2} \|h\|^2_2 \log\left(\frac{\|h\|^2_2}{\sum_{i=r+1}^n I^2_i}\right)
\ge \frac{2}{\log 2} \|h\|^2_2 \log\left(\frac{\|h\|^2_2}{K \alpha}\right)
$$
Recalling the definition of $\alpha$, it is now easy to choose the
error term $o_{\mu,\e}(1)$ appropriately, so that the last inequality
implies $\|h\|^2_2 \le \mu - 2\mu^2$.
\eprf

\subsection{The measure $\mu_p$}
Let $\mu_p$ be a product distribution, $\mu_p = 
\otimes_{k=1}^n \mu$, with $\mu(1) = p$, $\mu(0) = 1-p$. Assume $p \le
\frac12$.  
The best constant $C$ in the logarithmic Sobolev inequality in this
case is known \cite{d-sc} to be $C(p) = 
\frac{1-2p}{p(1-p)} \cdot \frac{1}{\log(1-p) - \log p}$.

\noindent {\bf Proof of a theorem of Friedgut and Kalai}
\thm \cite{friedgut-kalai}
Let $f$ be a boolean function with expectation $\mu$ on $\H^n$ endowed
with the measure $\mu_p$. Assume that $1 \gg p \ge n^{-o_n(1)}$. 
Then there is a variable with influence at least 
$\Omega\(\frac{\mu(1-\mu)}{p \log \frac 1p} \cdot
\frac{\log n}{n}\)$ on $f$.
\ethm
\prf
For a boolean function $f$ theorem~\ref{main} gives
$$
\label{bo_pr}
\sum_{i=1}^n I^2_i(f) \ge \frac{\sigma^2(f)}{4p^2(1-p)^2}  \exp
\left\{-\frac{1}{C(p)\sigma^2(f)} \sum_{i=1}^n I_i(f) \right\}   
$$
The expression on the right hand side is somewhat complicated. It
simplifies for $p \ll 1$, for which $C(p) \approx 
\frac{1}{p\log \frac 1p}$, and we get (ignoring negligible errors)
$$
\sum_{i=1}^n I^2_i(f) \ge \frac{\sigma^2(f)}{4p^2}  
\exp \left\{-\frac{p \log \frac 1p}{\sigma^2(f)} \sum_{i=1}^n I_i(f) \right\}
$$
Proceeding similarly to the proof of theorem~\ref{kkl}, we get that 
\beqn
\label{small_p}
\sum_{i=1}^n I^2_i(f) \ge \frac{\mu^2(1-\mu)^2}{p^2 \log^2 \frac 1p} \cdot
\frac{\log^2 n}{n}
\eeqn
In particular, there is a variable $i$ with influence at least
\beqn
\label{max_inf}
I_i \ge \frac{\mu(1-\mu)}{p \log \frac 1p} \cdot
\frac{\log n}{n}
\eeqn 
\eprf\\
We remark that this proof provides (\ref{max_inf}) with an explicit
constant $1$, and does not rely on the assumption $p \ge n^{-o(1)}$.

\section{Construction of 'isoperimetric' functions}
\label{Construction}
In this section we construct nonnegative functions $f_s$
on the cube endowed with the uniform measure, for which inequality
(\ref{ineq_main}) is almost tight. This directly implies that these functions
are 'isoperimetric' for inequalities (\ref{func_isop}),
(\ref{ineq_main_boolean}), and inequality (\ref{ineq_cnj_sobolev}) with
constant $c = 2$.

The functions $f_s$ were constructed in \cite{ns} (for a
different purpose). Here we repeat parts of this construction for completeness.

Let $s$ be an integer, $\sqrt{n} \ll s \ll n$. We first construct an
auxiliary function $k_s$. This function will be
{\it symmetric}, namely its value at a point will depend only on the
distance of the point from zero. Such
a function, of course, is fully defined by its values
$k_s(0),...,k_s(n)$ at distances $0...n$. Set $k_s(-1) = 0$ and
$k_s(0) = 1$, and define $k_s(r)$ for $1 \le r \le n$ so that the
relation $(n-2s) k_s(r) = r k_s(r-1) + (n-r) k_s(r+1)$ is satisfied
for $r = 0...n-1$. The univariate function $k_s(r)$ we have defined on the
integer points $r = 0...n$ coincides with a normalized Krawchouk polynomial
$K_s$ (see \cite{lev} for detailed information on Krawchouk polynomials). 

Krawchouk polynomials $\{K_s\}_{s=0}^n$ are a family of 
polynomials orthogonal with respect to a measure supported on
$0...n$. Hence their roots are simple and are located in the interval
$(0,n)$ \cite{szego}. Let $x_s$ be the first root of $K_s$. We now
define $f = f_s$ to be a symmetric function on
$\H^n$ defined by $f(x) = k_s(x)$ for points whose distance from zero
is at most $x_s$, and $f_s(x) = 0$ otherwise.  
 
We require an asymptotic estimate $x_s = \frac n2 - \sqrt{sn} +
o(\sqrt{sn})$ \cite{lev}. This means, in particular, that the support
of $f$ is small (of cardinality $e^{-s}$),
and therefore $\frac{\E^2 f}{\E f^2} \le 
\frac{\big |supp(f)\big |}{2^n}$ is small, so that $\sigma^2(f)$ can
be, for all practical reasons, replaced with $\E f^2$. 

For $x \in \H^n$, let $N(x) =
\sum_{y:~y \sim x} f(y)$. Then it is not hard to see (\cite{ns},
lemma~3.4) that $N(x) \ge (n-2s) f(x)$ and therefore ${\cal E}(f,f) = 
\E_x \sum_{y\sim x} (f(x) - f(y))^2 = 2n \E f^2 - 2 \<f,N\> \le
4s \E f^2$. 

Hence the right
hand side in (\ref{ineq_main}) can be estimated from below by $\E f^2
\cdot e^{-2sn}$. 

Now to the left hand side. The function $f$ is symmetric and (easy to
see, cf. also \cite{lev}) monotone. Therefore $\E |d_i| = I_i =
\hat{f}(\{i\})$ have the same value for all indices $1 \le i\le n$. Let $I$
denote this common value. We want to upper bound $I$. Let $f(x)$
denote the value of $f$ in points at distance $x$ from zero. Let $m =
\lfloor x_s \rfloor$. Then
$$
n I = \sum_{i=1}^n \hat{f}\(\{i\}\) =
\frac{1}{2^{2n}} \sum_{x=0}^n {n\choose x} (n-2x) f(x) = \frac{1}{2^{2n}}
\sum_{x=0}^m {n\choose x} (n-2x) f(x) \le \frac{n}{2^{2n}} \sum_{x=0}^m
    {n\choose x} f(x)
$$    
Since $\E f^2 = \frac{1}{2^n} \sum_{x=0}^n {n\choose x} f^2(x)$, we
have $f(x) \le \sqrt{\frac{2^n \E f^2}{{n \choose x}}}$, and so
$$
\frac{1}{2^n} \sum_{x=0}^m {n\choose x} f(x) \le \sqrt{2^n \E f^2} 
\sum_{x=0}^m \sqrt{ {n\choose x}} \le n \sqrt{{n\choose m} 2^n \E f^2}
$$

Therefore $I \le n \sqrt{\frac{{n\choose m} \E
f^2}{2^n}}$, and the left hand side of (\ref{ineq_main}) is
$n I^2 \le n^2 \frac{{n\choose m}}{2^n} \E f^2$. Now observe
\cite{lint} that ${n\choose m} \le 2^{nH\(m/n\)} \le 2^n
e^{-(1-o_n(1))2sn}$. 

Hence the left hand side in (\ref{ineq_main}) can be estimated from
below by $\E f^2 \cdot n^2 e^{-(1-o_n(1))2sn}$.
The estimates for both sides are sufficiently close to show the constant
$1/2$ in the exponent on 
the right hand side of (\ref{ineq_main}) to be best possible.

\section{Appendix A - an inductive proof of an isoperimetric
  inequality}
\label{Appendix-A}
In this section we give an inductive proof of the inequality (\ref{func_isop})
$$
{\cal E}(f,f) \ge 2 \cdot \E f^2
\log \frac{\E f^2}{\E^2 f}
$$
for a real nonnegative function $f~:\H^n \rarrow \R$.
By homogeneity, we may and will assume $\E f = 1$.

The proof is by induction on the dimension $n$.

For $n = 1$, let
$f(0) = a$, and $f(1) = 2-a$. Then $\E f^2 = \frac12 \(a^2 + (2-a)^2\)
= 1 + (1-a)^2$; and $\E_x \sum_{y\sim x} (f(x) - f(y))^2 = (2-2a)^2
= 4(1-a)^2$. Let $x = (1 - a)^2$. It remains to verify
that $2x \ge (1+x) \log (1+x)$ for $0\le x\le 1$, which is easily
seen to be true. In fact a stronger inequality $2x \ge (1+x)
\log_2(1+x)$ is also valid in this interval,  since the right hand
side is a convex function which is $0$ at zero and $2$ at one.

Assume the inequality to hold for $n-1$.
Let $f_0$ and $f_1$ be the
restrictions of $f$ to $(n-1)$-dimensional half-cubes determined
by value of the $n$-th coordinate. Let $\mu_i$ be the expectations of $f_i$,
and $v_i$ the second moments of $f_i$ for $i=0,1$.

Then 
\beqn
\label{induct}
{\cal E}(f,f) = 
\frac 12 \cdot \left({\cal E}\(f_0,f_0\) + {\cal E}\(f_1,f_1\)\right)
+ \|f_0 - f_1\|^2_2. 
\eeqn
Note that the expectations and the distance in this formula are
computed on $(n-1)$-dimensional cubes.
 
By the induction hypothesis we can lower bound the first summand by 
$$
v_0 \log \frac{v_0}{\mu^2_0} + v_1 \log \frac{v_1}{\mu^2_1}. 
$$ 
For the second summand we need a simple lemma.
\lem
Let $f_0$ and $f_1$ be two functions with expectations $\mu_0, \mu_1$
and variances $\sigma^2_0 = v_0 - \mu^2_0$, $\sigma^2_1 = v_1 -
\mu^2_1$. Then 
$$
\|f_0 - f_1\|^2 \ge (\sigma_0 - \sigma_1)^2 + (\mu_0 - \mu_1)^2.
$$
\elem
\prf 
Let $g_i = f_i - \mu_i$, $i = 0,1$. Then $\E g_i = 0$ and therefore 
$$
\|f_0 - f_1\|^2 = \left<g_0 - g_1 ,g_0 - g_1\right> + (\mu_0-\mu_1)^2
=  \|g_0 - g_1\|^2 + (\mu_0-\mu_1)^2 \ge (\sigma_0 - \sigma_1)^2  +
(\mu_0-\mu_1)^2. 
$$
\eprf

Going back, and substituting in (\ref{induct}), 
$$
{\cal E}(f,f) \ge 
v_0 \log\(\frac{v_0}{\mu^2_0}\) + v_1 \log\(\frac{v_1}{\mu^2_1}\) +
\left[v_0 + v_1 - 2\sqrt{v_0 - \mu^2_0}\sqrt{v_1 - \mu^2_1} - 2 \mu_0
  \mu_1\right]. 
$$
So it suffices to show that under the assumptions
\begin{enumerate}
\item
$\mu_0, \mu_1, v_0, v_1 \ge 0$,
\item
$\frac{\mu_0 + \mu_1}{2} = 1$, and
\item
$\frac{v_0 + v_1}{2} = v := \E f^2$
\end{enumerate}
holds
\beqn
\label{main_func_isop_proof}
v_0 \log\(\frac{v_0}{\mu^2_0}\) + v_1 \log\(\frac{v_1}{\mu^2_1}\) +
\left[v_0 + v_1 - 2\sqrt{v_0 - \mu^2_0}\sqrt{v_1 - \mu^2_1} - 2 \mu_0
  \mu_1\right] \ge 2v \log v.
\eeqn
The next few steps swap variables to simplify this expression.

Take $t = \frac{\mu_0-\mu_1}{2}$. Then $\mu_0 =
1 + t$ and $\mu_1 = 1 - t$. Similarly take $v_0 = v(1+y)$ and $v_1 =
v(1-y)$. Note that $-1 \le t, y\le 1$. 

Substituting in (\ref{main}), and dividing out by $2v$ it needs to be seen that
$$
\frac{1+y}{2} \log\left(\frac{1+y}{(1+t)^2}\cdot v\right) + 
\frac{1-y}{2} \log\left(\frac{1-y}{(1-t)^2}\cdot v\right) + 1 
\ge  
$$
$$
\log v + \frac{1}{v} \cdot \left[\sqrt{v(1+y) - (1+t)^2}\sqrt{v(1-y) -
    (1-t)^2} + (1-t^2)\right],
$$
or
$$
\frac{1+y}{2} \log\left(\frac{1+y}{(1+t)^2}\right) + 
\frac{1-y}{2} \log\left(\frac{1-y}{(1-t)^2}\right) + 1 
\ge 
$$
$$
\frac{1}{v} \cdot \left[\sqrt{v(1+y) - (1+t)^2}\sqrt{v(1-y) -
    (1-t)^2} + (1-t^2)\right].
$$
We first take on the right hand side and show it to be at most $\sqrt{1-y^2}$. 
Indeed, it suffices to show
$$
\sqrt{v^2 \(1-y^2\) - v(1+y)(1-t)^2 - v(1-y)(1+t)^2 + \(1-t^2\)^2} \le
v\sqrt{1-y^2}  - \(1-t^2\). 
$$ 
Note that the right hand side is nonnegative, since going back to the
definitions of $t$ and $y$, this is $\sqrt{v_0 v_1} - \mu_0 \mu_1$. 
Squaring both expressions, and rearranging, we get to
$$
(1+y)(1-t)^2 + (1-y)(1+t)^2 \ge 2 \sqrt{1-y^2}\(1-t^2\).
$$
This inequality is a special case of the Arithmetic-Geometric inequality.

Now to the left hand side. Let $H(x) = -x \log x -
(1-x) \log (1-x)$ be the (natural) entropy function. For $0 \le p, q\le
1$ let $D(p||q) = p \log \frac pq + (1-p) \log
\frac{1-p}{1-q}$ denote the divergence between two-point distributions
$(p,1-p)$ and $(q,1-q)$. It is well-known (and is a simple consequence of
the concavity of logarithm) that divergence is nonnegative. Now,
$$
1 + \left[\frac{1+y}{2} \log\left(\frac{(1+y)^2}{(1+t)^2}\right) + 
\frac{1-y}{2} \log\left(\frac{(1-y)^2}{(1-t)^2}\right)\right] -
\left[\frac{1+y}{2} \log(1+y) + \frac{1-y}{2} \log(1-y)\right] = 
$$
$$
1 + 2 \left[\frac{1+y}{2} \log\left(\frac{(1+y)}{(1+t)}\right) + 
\frac{1-y}{2} \log\left(\frac{(1-y)}{(1-t)}\right)\right] +
H\(\frac{1-y}{2}\) - \log 2 = 
$$
$$
2 D\(\frac{1-y}{2} \Big | \Big |\frac{1-t}{2}\) + 
H\(\frac{1-y}{2}\) + (1 - \log 2) \ge 
H\(\frac{1-y}{2}\) + (1 - \log 2).
$$
Therefore we need to show
$$
1 - \sqrt{1-y^2} \ge \log 2 - H\(\frac{1-y}{2}\)
$$
for all $-1\le y\le 1$.

We will need two well-known facts: the function $\phi(t) = \frac12 -
\sqrt{t(1-t)}$ is an involution on $[0,\frac12]$; and the function
$R(x) = H(\phi(x))$ is convex on $[0,\frac12]$.

Since $R(0) = \log 2$ and $R'(0) = -2$, by convexity $R(x)
\ge \log 2 - 2 x$ for $x \in [0,\frac12]$. Rearranging and taking $x =
\phi(z)$,  
$$
2\phi(z) \ge \log 2 -H(z)
$$
Substituting $y = 1 - 2z$
$$
1 - \sqrt{1-y^2} \ge \log 2 - H\(\frac{1-y}{2}\)
$$
and we are done.

\section{Appendix B - an inductive proof of a logarithmic Sobolev
  inequality}
\label{Appendix-B}
In this section we give an inductive proof of the inequality
$$
{\cal E}(f,f) \ge 2 \cdot \sum_{i=1}^n Ent(d^2_i)
$$
for a real function $f~:\H \rarrow \R$, with $d_i$ the difference
functions of $f$ (cf. section~\ref{Main}).   

Observe that the right hand side might depend on the ordering of coordinates. 

The proof is by induction on $n$. For $n=1$, $d_1 = f - \E f$, and
therefore $d^2_1$ is a constant function with zero entropy. The claim
follows. Assume the claim for $n-1$,
and consider it for $n$. 

Let the influence $I_i(f)$ of the
$i$-th bit on a real-valued function $f$ to be given by $\E_x (f(x) -
f(x \oplus e_i))^2$. Thus ${\cal E}(f,f) = \sum_{i=1}^n I_i$.  

Let $f_0$, $f_1$ be the restrictions of $f$ to subcubes defined by the
value of the $n$-th coordinate. We write $f \leftrightarrow (f_0,f_1)$. These
are functions on $n-1$ variables. Let 
their influences, their conditional expectations, and their difference
functions with respect to the 
natural ordering $1...n-1$ of the coordinates be denoted by $I_{i,0}$,
$I_{i,1}$, $f_{i,0}$, $f_{i,1}$, $d_{i,0}$, $d_{i,1}$
correspondingly. Then, by the induction hypothesis 
$$
\sum_{i=1}^n I_i = \frac12 \(\sum_{i=1}^{n-1} I_{i,0} +
\sum_{i=1}^{n-1} I_{i,1} \) + I_n \ge
\sum_{i=1}^{n-1} Ent\(d^2_{i,0}\) + \sum_{i=1}^{n-1}
Ent\(d^2_{i,1}\) + I_n
$$
Consider now a slightly different ordering of the coordinates for $f$, which
is $n,1,2...,n-1$. Let $f_i$ and $d_i$ be the conditional expectations and the
difference functions in this new ordering. We will show
$$
\sum_{i=1}^{n-1} Ent\(d^2_{i,0}\) + \sum_{i=1}^{n-1}
Ent\(d^2_{i,1}\) + I_n \ge 2 \cdot \sum_{i=1}^n Ent\(d^2_i\)
$$
This will prove the inequality for the ordering $n,1,2...,n-1$ of the
coordinates.\footnote{Alternatively, we could have insisted on the
  'natural' order of 
  the coordinates for $f$, and changed the order of coordinates for
  $f_0$, $f_1$.}

Observe $f_i \lrarrow (f_{i-1,0},~f_{i-1,1})$ for $i=2...n$, and
similarly for difference functions. Note also that $d^2_1$ is
constant, and therefore has zero entropy.
\lem
Let $k\ge 0$, $k \lrarrow (g,h)$. Then 
$$
Ent(k) = \frac12 (Ent(g) + Ent(h)) + \frac12\(\E g \log \E g + \E h \log
\E h - (\E g + \E h) \log \frac{\E g + \E h}{2}\)
$$
\elem
\prf
\eprf\\
Therefore we need to show
$$
\frac 12 \cdot I_n \ge \sum_{i=2}^n \(Ent\(d^2_i\) -
\frac12\(Ent\(d^2_{i-1,0}\) + Ent\(d^2_{i-1,1}\)\)\) = 
$$
$$
\frac12 \sum_{i=1}^{n-1} \( \a_i \log \a_i + \b_i \log \b_i - \(\a_i + \b_i\)
\log \frac{\a_i + \b_i}{2}\)
$$
Wor this purpose we need information on the joint behaviour of the sequences
$\a_i$ and $\b_i$.  

Applying the Cauchy-Schwarz inequality twice,
$$
\(\a_i - \b_i\)^2 = \E^2\(d_{i,0} - d_{i,1}\)\(d_{i,0} +
d_{i,1}\)  \le \E^2 |d_{i,0} - d_{i,1}| |d_{i,0} + d_{i,1}| \le
$$
$$
\E\(d_{i,0} - d_{i,1}\)^2 \E\(d_{i,0} + d_{i,1}\)^2 \le
\(\a_i + \b_i + 2\sqrt{\a_i \b_i}\) \cdot \E\(d_{i,0} - d_{i,1}\)^2
$$
Observe that $d_{i,0} - d_{i,1}$ is a difference sequence for $f_0 -
f_1$ and therefore $\sum_{i=1}^{n-1} \E\(d_{i,0} - d_{i,1}\)^2 \le
(f_0 - f_1)^2 = I_n$. 

Take $a_i := \a_i$, $b_i := \b_i$, and consider an optimization
problem
$$
\mbox{Maximize}~~\sum_{i=1}^{n-1} \(a_i \log a_i + b_i \log b_i - (a_i
+ b_i) \log \frac{a_i + b_i}{2}\)
$$
$$
\mbox{Given}~~ \sum_{i=1}^{n-1} \frac{(a_i - b_i)^2}{a_i + b_i +
  2\sqrt{a_i b_i}} \le I_n,~~~a_i, b_i \ge 0
$$ 
Let $m_i = min\{a_i, b_i\}$ and $c_i = |a_i - b_i|$. Then an
equivalent formulation is
$$
\mbox{Maximize}~~\sum_{i=1}^{n-1} \((m_i + c_i) \log (m_i + c_i) + m_i
\log m_i - (2m_i + c_i) \log 
\frac{2m_i + c_i}{2}\)
$$
$$
\mbox{Given}~~ \sum_{i=1}^{n-1} \frac{c^2_i}{2m_i + c_i +
  2\sqrt{m_i(m_i+ c_i)}} \le I_n,~~~m_i, c_i \ge 0
$$
Assume $c_i > 0$ for all $i$ since removing coordinates with $c_i = 0$
does not effect neither the target function nor the
constraint. Therefore we are allowed to consider $r_i =
\frac{m_i}{c_i}$, leading to the following formulation
$$
\mbox{Maximize}~~\sum_{i=1}^{n-1} c_i \cdot \((1 + r_i) \log (1+r_i) +
r_i \log r_i - (1 + 2 r_i) \log \frac{1 + 2 r_i}{2} \)
$$
$$
\mbox{Given}~~ \sum_{i=1}^{n-1} \frac{c_i}{1 + 2 r_i +
  2\sqrt{r_i(1+ r_i)}} \le I_n,~~~r_i, c_i \ge 0
$$
The following technical claim completes the analysis.
\lem
For any $r\ge 0$ holds
$$
(1 + r) \log (1+r) + r \log r - (1 + 2 r) \log \frac{1 + 2 r}{2} \le 
\frac{1}{1 + 2 r + 2\sqrt{r(1+ r)}}
$$
In addition 
$$
\(1 + 2 r + 2\sqrt{r(1+ r)}\) \cdot \((1 + r) \log (1+r) + r \log r - (1
+ 2 r) \log \frac{1 + 2 r}{2}\) \rarrow_{r\rarrow \infty} 1
$$
\elem
Therefore the supremum of the above maximization problem is bounded by $I_n$
(and it can easily be seen that it actually equals $I_n$), completing
the proof.\\
\prf (Of the lemma)\\
Let $g(r) = 1 + 2 r + 2\sqrt{r(1+ r)} = \(\sqrt{r} + \sqrt{1+r}\)^2$,
and $h(r) = (1 + r) \log (1+r) + r \log r - (1
+ 2 r) \log \frac{1 + 2 r}{2}$. We want to show $(gh)(r) \le 1$, for all
$r\ge 0$. At zero, $(gh)(0) = \log 2 < 1$, at infinity, $g(r) \sim 4r$
and $h(r) \sim \frac{1}{4r}$, and thus $(gh)(r) \rarrow_{r \rarrow
  \infty} 1$, proving the second part of the lemma. Thus it is
sufficient to show $gh$ is increasing, or $\frac{g'h}{g} \ge -h'$.

Computing, $\frac{g'h}{g} = \frac{h}{\sqrt{r(1+r)}}$, and $-h' =
\log\frac{(1+2r)^2}{4r(1+r)}$. It remains to show
$
h(r) \ge \sqrt{r(1+r)} \cdot \log\frac{(1+2r)^2}{4r(1+r)}
$
Rewriting $h(r)$ as $\log \frac{2+2r}{1+2r} - r
\log\frac{(1+2r)^2}{4r(1+r)}$, this is the same as 
$$
\log \frac{2+2r}{1+2r} \ge \(\sqrt{r(1+r)} + r\) \cdot
\log \frac{(1+2r)^2}{4r(1+r)} = \frac{r}{\sqrt{r(1+r)} - r}\cdot
\log \frac{(1+2r)^2}{4r(1+r)}
$$
Or
$$
\(\sqrt\frac{1+r}{r} - 1\) \cdot \log \frac{2+2r}{1+2r} \ge
\log\frac{(1+2r)^2}{4r(1+r)}.
$$
Let $t = \sqrt\frac{1+r}{r}$. Then $t \in (1,\infty)$. Rewriting in
terms of $t$, we want to have 
$$
(t-1) \log \frac{2t^2}{t^2 + 1} \ge 2 \log \frac{t^2 + 1}{2t}
$$
This holds at one. Comparing the derivatives, it suffices to show 
$$
\log \frac{2t^2}{t^2 + 1} \ge \frac{2t-2}{t^2 + 1}.
$$
Once again, this holds at one. Comparing the derivatives for the final
time, one has to show
$$
\frac 1t \ge \frac{-t^2 + 2t + 1}{t^2 + 1},
$$
or $t^3 + 1 \ge t^2 + t$, which is immediate for $t\ge 1$.
\eprf

\section{Acknowledgements}
We are grateful to Ehud Friedgut for his suggestions which led
to a significant simplification of the proof of theorem~\ref{main}. 
We also thank Gil Kalai, Nati Linial, and Amites Sarkar for many valuable
remarks.

\end{document}